%
\magnification=\magstep1
\vsize 23 truecm   
\hsize 15.5 truecm 
\baselineskip=13.5 pt 
\parskip=3pt 
\overfullrule=0pt


%
\def\BBF#1{\expandafter\edef\csname #1#1\endcsname{%
   {\mathord{\bf #1}}}}
\BBF m
\BBF p
\BBF C
\BBF K
\BBF N
\BBF P
\BBF Q
\BBF R
\BBF T
\BBF V
\BBF Z

%
\def\CAL#1{\expandafter\edef\csname #1\endcsname{%
   {\mathord{\cal #1}}}}
\CAL A
\CAL B
\CAL C
\CAL D
\CAL E
\CAL F
\CAL G
\CAL H
\CAL{IH}
\CAL J
\CAL K
\CAL L
\CAL P

%
\newread\amssymdef
\openin\amssymdef=amssym.def
\ifeof\amssymdef\else\input amssym.def
\def\mm{\frak m}

\fi
\closein\amssymdef

\newread\amssymtex
\openin\amssymtex=amssym.tex
\ifeof\amssymtex\def\precneqq{%
\mathrel{\mathop\prec\limits_{\scriptscriptstyle\ne}}}
\else\input amssym.tex
\fi
\closein\amssymtex

%
\def\:{\colon}

\def\b{}

\def\quer{\overline}

%
\def\longto{\longrightarrow}

\def\onto{\to\kern-.8em\to}
\def\longonto{\longto\kern-.8em\longto}
\def\longmapsto{\mapstochar\longto}

%
\def\epsilon{\varepsilon}
\def\phi{\varphi}
\def\rho{\varrho}
\def\theta{\vartheta}

%
\def\widecong{\;\cong\;}
\def\wideto{\;\to\;}

%
\def\qed{\hfill$\scriptscriptstyle\rfloor$}
\def\sqr#1#2{{\vbox{\hrule height.#2pt
         \hbox{\vrule width.#2pt height#1pt \kern#1pt
               \vrule width.#2pt}
         \hrule height.#2pt}}}

\def\qed{\ifmmode\hbox{\sqr35}
     \else\hfill\sqr35\fi}

%
\def\ibull#1\par{{\parindent10pt\item{$\bullet$}#1\par}}
\def\iibull#1\par{{\parindent10pt\itemitem{$\bullet$}#1\par}}
\def\iiitem{\par\indent\indent\hangindent3\parindent\textindent}
\def\iiibull#1\par{{\parindent10pt\iiitem{$\bullet$}#1\par}}
\def\istern#1\par{{\parindent10pt\item{$*$}#1\par}}
\def\istrich#1\par{{\parindent10pt\item{--}#1\par}}

\def\inditem#1#2\par{{\parindent10pt\item{#1}#2\par}}

%
\def\MOPnl#1{\expandafter\edef\csname #1\endcsname{%
   {\mathop{\rm #1}\nolimits}}}
\MOPnl{cld}
\MOPnl{codim}
\MOPnl{coker}
\MOPnl{Coker}
\MOPnl{Ext}
\MOPnl{Hom}
\MOPnl{Ker}
\MOPnl{id}
\MOPnl{im}
\MOPnl{ker}
\MOPnl{lin}
\MOPnl{max}
\MOPnl{odd}
\MOPnl{or}
\MOPnl{pr}
\MOPnl{pt}
\MOPnl{rg}
\MOPnl{rk}
\MOPnl{Sp}
\MOPnl{Spec}
\MOPnl{st}
\MOPnl{supp}
\MOPnl{Tor}

%
\font\XIIbfM=cmbx12 scaled 1200
\font\XIIbf=cmbx12

\font\XIIrm=cmr12

\font\bfit=cmbxti10
\font\sc=cmcsc10

%
\def\diagram{\def\normalbaselines{\baselineskip20pt\lineskip3pt
\lineskiplimit3pt}\matrix}
\def\mapright#1{\smash{\mathop{\hbox to 35pt{\rightarrowfill}}\limits^{#1}}}
\def\mapleft#1{\smash{\mathop{\hbox to 35pt{\leftarrowfill}}\limits^{#1}}}
\def\mapdown#1{\Big\downarrow\rlap{$\vcenter{\hbox{$\scriptstyle#1$}}$}}

%

\everydisplay{\abovedisplayskip 4.5pt plus 3 pt minus 1.5 pt
\belowdisplayskip\abovedisplayskip
\belowdisplayshortskip=3pt plus 3pt minus 1pt
\thickmuskip=1.5\thickmuskip plus 5.0mu
}

%
%
\font\sevenrm=cmr7     \font\sixrm=cmr6
\font\seveni=cmmi7     \font\sixi=cmmi6
\font\sevensy=cmsy7    \font\sixsy=cmsy6
\font\sevenbf=cmbx7    \font\sixbf=cmbx6
\font\seventt=cmtt8 scaled 875
\font\fivei=cmmi5
\font\sevenit=cmti7
\font\sevensl=cmsl8 scaled 875

\font\fivebf=cmbx5
\def\kleinst{\def\rm{\fam0\sevenrm}
  \textfont0=\sevenrm \scriptfont0=\sixrm \scriptscriptfont0=\fiverm
  \textfont1=\seveni   \scriptfont1=\sixi  \scriptscriptfont1=\fivei
  \textfont2=\sevensy   \scriptfont2=\sixsy \scriptscriptfont2=\fivesy
  \textfont\itfam=\sevenit \def\it{\fam\itfam\sevenit}%
  \textfont\slfam=\sevensl \def\sl{\fam\slfam\sevensl}%
 \textfont\ttfam=\seventt \def\tt{\fam\ttfam\seventt}%
  \textfont\bffam=\sevenbf \scriptfont\bffam=\sixbf
   \scriptscriptfont\bffam=\fivebf \def\bf{\fam\bffam\sevenbf}%
  \normalbaselineskip=7pt
  \setbox\strutbox=\hbox{\vrule height6pt depth2pt width0pt}%
  \let\sc=\sixrm \let\big=\sevenbig \normalbaselines\rm}
\def\sevenbig#1{{\hbox{$\textfont0=\eightrm\textfont2=\eightsy
    \left#1\vbox to6.5pt{}\right.\n@space$}}}

%
\font\eightrm=cmr8     
\font\eighti=cmmi8     
\font\eightsy=cmsy8    
\font\eightbf=cmbx8    
\font\eighttt=cmtt8    
\font\eightit=cmti8
\font\eightsl=cmsl8
\def\kleiner{\def\rm{\fam0\eightrm}
  \textfont0=\eightrm \scriptfont0=\sixrm \scriptscriptfont0=\fiverm
  \textfont1=\eighti     \scriptfont1=\sixi \scriptscriptfont1=\fivei
  \textfont2=\eightsy     \scriptfont2=\sixsy \scriptscriptfont2=\fivesy
  \textfont\itfam=\eightit \def\it{\fam\itfam\eightit}%
  \textfont\slfam=\eightsl \def\sl{\fam\slfam\eightsl}%
 \textfont\ttfam=\eighttt \def\tt{\fam\ttfam\eighttt}%
  \textfont\bffam=\eightbf \scriptfont\bffam=\sixbf
   \scriptscriptfont\bffam=\fivebf \def\bf{\fam\bffam\eightbf}%
  \normalbaselineskip=8pt
  \setbox\strutbox=\hbox{\vrule height6pt depth2pt width0pt}%
  \let\sc=\sixrm \let\big=\eightbig \normalbaselines\rm}

\def\eightbig#1{{\hbox{$\textfont0=\ninerm\textfont2=\ninesy
    \left#1\vbox to6.5pt{}\right.\n@space$}}}
\def\ninebig#1{{\hbox{$\textfont0=\tenrm\textfont2=\tensy
    \left#1\vbox to7.25pt{}\right.\n@space$}}}

%
%
\font\ninerm=cmr9
\font\ninei=cmmi9
\font\ninesy=cmsy9
\font\ninebf=cmbx9
\font\ninett=cmtt9
\font\nineit=cmti9
\font\ninesl=cmsl9
\def\klein{\def\rm{\fam0\ninerm}
  \textfont0=\ninerm \scriptfont0=\sixrm \scriptscriptfont0=\fiverm
  \textfont1=\ninei     \scriptfont1=\sixi \scriptscriptfont1=\fivei
  \textfont2=\ninesy     \scriptfont2=\sixsy \scriptscriptfont2=\fivesy
  \textfont\itfam=\nineit \def\it{\fam\itfam\nineit}%
  \textfont\slfam=\ninesl \def\sl{\fam\slfam\ninesl}%
  \textfont\ttfam=\ninett \def\tt{\fam\ttfam\ninett}%
  \textfont\bffam=\ninebf \scriptfont\bffam=\sixbf
  \scriptscriptfont\bffam=\fivebf \def\bf{\fam\bffam\ninebf}%
  \normalbaselineskip=10pt
  \setbox\strutbox=\hbox{\vrule height7pt depth3pt width0pt}%
  \let\sc=\sevenrm \let\big=\ninebig \normalbaselines\rm}

%
%
    \newcount\annotation
     \annotation=0
      \long\def\Fussnote#1{{\baselineskip=9pt
           \setbox\strutbox=\hbox{\vrule height 7pt depth 2pt width 0pt}
         \klein\global\advance\annotation by 1
         \footnote{$^{\the\annotation)}$}{#1}}}
	 
%
%
\def\makeatletter{\catcode`\@=11\relax}
\def\makeatother{\catcode`\@=12\relax}
\makeatletter

\newdimen\@tempdima
\newbox\@tempboxa

\def\@height{height}
\def\@depth{depth}
\def\@width{width}

\newdimen\fboxrule
\newdimen\fboxsep

\fboxrule = .4pt
\fboxsep = 2pt

\long\def\fbox#1{\leavevmode\setbox\@tempboxa\hbox{#1}\@tempdima\fboxrule
   \advance\@tempdima \fboxsep \advance\@tempdima \dp\@tempboxa
   \hbox{\lower \@tempdima\hbox
  {\vbox{\hrule \@height \fboxrule
          \hbox{\vrule \@width \fboxrule \hskip\fboxsep
          \vbox{\vskip\fboxsep \box\@tempboxa\vskip\fboxsep}\hskip
                 \fboxsep\vrule \@width \fboxrule}%
                 \hrule \@height \fboxrule}}}}

\makeatother

\def\LEB{}\long\def\TOT#1\LEB{}
\headline{\ifnum\pageno=1\else\klein\hfill 
Combinatorial Duality and Intersection Product 
\hfill \folio\ \fi}
\footline{}

\bigskip\noindent


\centerline{\XIIbfM Combinatorial Duality and Intersection Product:}
\vskip 6 pt
\centerline{\XIIbfM A Direct Approach}
\vskip 12 pt
\centerline{\XIIrm Gottfried Barthel, Jean-Paul Brasselet,}
\vskip 3 pt
\centerline{\XIIrm Karl-Heinz Fieseler, Ludger Kaup}

\vskip 0.8 true cm
{\klein\narrower\noindent
{\bf Abstract:} 
The proof of the combinatorial Hard Lefschetz Theorem for the 
``virtual'' intersection cohomology of a not necessarily rational 
polytopal fan that has been presented in [Ka] completely establishes 
Stanley's conjectures for the generalized $h$-vector of an arbitrary 
polytope. The main ingredients, namely, Poincar\'e Duality and the 
Hard Lefschetz Theorem, both rely on the intersection product. In the 
constructions of [BBFK${}_2$] and [BreLu$_{1}$], there remained an apparent 
ambiguity. The recent solution of this problem in [BreLu$_2$] uses the 
formalism of derived categories. The present article gives a 
straightforward approach to combinatorial duality and a natural 
intersection product, completely within the framework of elementary 
sheaf theory and commutative algebra, thus avoiding derived categories. 
\par}

\bigskip
\bigskip
\centerline{\XIIbf Introduction}

\medskip
\noindent
In [St], R.~Stanley introduced the generalized $h$-vector for arbitrary 
polytopes. For rational polytopes, this new combinatorial invariant 
agrees with the vector of even (middle perversity) intersection 
cohomology Betti numbers of a projective toric variety associated with 
the polytope and then, they enjoy the same properties. Stanley proved that the
Dehn-Sommerville equations (i.e.,  Poincar\'e duality) remain valid in the
general case, and he conjectured  that non-negativity and unimodality up to the
middle dimension also  should continue to hold. In the rational case, the
unimodality property  follows from the ``Hard Lefschetz Theorem'' for the
intersection  cohomology of a projective variety.

This conjecture motivated the search for a purely combinatorial 
approach to the intersection cohomology of toric varieties. 
Such an approach has been developed independently in [BBFK$_{2}$] and 
by P.~Bressler and V.~Lunts in [BreLu$_{1}$]. The basic idea in both 
articles is to view a (not necessarily rational) fan as a finite 
topological space, endowed with the topology given by the subfans as 
non-trivial open sets, and to study the properties of a certain sheaf 
on that ``fan space'' that agrees with the equivariant intersection 
cohomology sheaf for the associated toric variety in the rational case. 
This approach then yields a ``virtual'' intersection cohomology theory 
for the class of ``quasi-convex'' fans that includes all complete and 
hence, in particular, all polytopal fans. In [BBFK$_{2}$], the working 
principle was to present everything on a fairly elementary level, using 
geometry and commutative algebra only and avoiding the use of derived 
categories.

At the time when these articles were written, a purely combinatorial 
version of the Hard Lefschetz Theorem (HLT), as stated in the third 
section, was still lacking. This was the only missing piece to prove that 
the vector of even ``virtual'' intersection cohomology Betti numbers of a 
polytopal fan agrees with the generalized $h$-vector of the polytope, and 
thus, to fully establish Stanley's conjecture. As another problem, in the 
construction of the 
intersection product on the virtual equivariant intersection cohomology
sheaf, apparently non-canonical choices entered.

In the meantime, a proof of the combinatorial Hard Lefschetz Theorem 
has been presented by K.~Karu in [Ka]. The proof heavily relies on the study 
of the intersection product, since what actually is shown are the 
Hodge-Riemann bilinear relations (HRR) for the ``primitive'' (virtual) 
intersection cohomology, which imply HLT as an easy consequence. The 
apparent ambiguity in the definition of the intersection product, 
however, makes the argumentation quite involved, since one has to 
carefully keep track of the choices made. A first simplified version 
has recently been presented by Bressler and Lunts in [BreLu$_2$], 
using the framework of derived categories. In parti\-cular, they verify 
by a detailed analysis that none of the possible choices affects the 
definition of the pairing.

Our goal is to go one step further, namely, to give a short, direct,  and
elementary approach to duality and the intersection product in the 
``geometrical'' spirit of [BBFK$_{2}$], following ideas of [Bri], the  only
prerequisites being sheaf theory and commutative algebra. While in
[BreLu${}_1$,2] the duality functor is apriori only defined as an endofunctor of
a big derived category containing the "pure sheaves" as invariant subcategory, we
construct the dual of a pure sheaf directly as a pure sheaf, avoiding the above
detour. The crucial step is the definition of the restriction homomorphisms from 
a cone to a facet; the image of a section can be looked at as a kind  of residue
along the facet. The structure of the proof of Poincar\'e duality and the
naturality of the intersection product is the classical one, as in [BreLu${}_2$],
with the single steps easily accessible. For the convenience of the reader we
give here a complete presentation, referring always to the corresponding
statements in [BreLu${}_1$,2]. An intersection product corresponds to a 
sheaf homomorphism $\theta \: \E\b \to \D\E\b$ of degree zero from the 
(graded) equivariant intersection cohomology sheaf $\E\b$ to its dual, 
and we check that $\E\b$ is self-dual in a natural way.\par

Our main result, cf also [BreLu${}_2$] is stated below in such a way that it fits
into the  inductive proof of the Hard Lefschetz Theorem as given in [Ka]: 
Assuming HLT for polytopal fans 
in dimension $d<n$, the ``Poincar\'e Duality Theorem'' yields a natural
intersection product on every fan in  dimension $n$. In [Ka] it is shown
that HRR for simplicial fans  in any dimension -- which is valid by [Mc] --
together with HRR for  arbitrary fans in dimensions $d<n$ imply HRR and thus, HLT 
in dimension~$n$. 
In that induction step, it is most useful to work 
with a canonical pairing.
  
To state our result, we use this notation, explained more systematically 
in section~0: Let~$\Delta$ be a quasi-convex fan in a vector space~$V$ 
of real dimension~$n$ with a fixed volume form, and $\partial\Delta$ 
its boundary fan. The global sections of the equivariant intersection 
cohomology sheaf $\E\b$ on~$\Delta$ and on $(\Delta, \partial\Delta)$, 
respectively, are (graded) modules over the (graded) symmetric algebra 
$A\b := S\b(V^*)$ of polynomial functions on~$V$.

\medskip
\noindent
{\bf Poincar\'e Duality Theorem.} [BreLu${}_2$, 3.16] {\it In the above
situation, let  us assume that the Hard Lefschetz Theorem holds in all dimensions 
below~$n$. Then there is a natural intersection product
$$
  \E\b (\Delta) \times \E\b (\Delta, \partial\Delta)
  \;\;\longto\;\; A\b[-2n]\leqno{\rm(PD)}
$$
giving rise to a dual pairing of finitely generated free $A\b$-modules.}

\smallskip
For the following supplement, let $\hat\E\b$ be the equivariant 
intersection cohomology sheaf of a refinement $\hat \Delta$ of $\Delta$ 
with refinement map $\pi: \hat \Delta \to \Delta$. The result is 
essential for a simplified proof of the Hard Lefschetz Theorem, see also 
[BreLu$_2$]:

\medskip
\noindent
{\bf Compatibility Theorem.} {\it Let $\E\b \to \pi_* (\hat \E\b)$ be a 
homomorphism of graded sheaves extending the identity $\E\b_{o} = \RR = 
\pi_*(\hat \E\b)_{o}$ at the zero cone~$o$. Then the ``global'' 
intersection products are compatible, i.e., the following diagram is 
commutative:
$$
\matrix{
\E\b (\Delta) \times \E\b (\Delta, \partial\Delta)\kern-30pt &&
  \longto && 
  \kern-30pt\hat \E\b (\hat\Delta) \times 
  \hat \E\b (\hat\Delta,\partial \hat \Delta)\cr
  \noalign{\vskip6pt}
  & \searrow && \swarrow & \cr 
  \noalign{\vskip6pt}
  & & A\b[-2n]\;. &  &\cr} 
$$
}

\medskip
The present article is a complete version of the results announced 
in~[Fi]. 
-- The authors gratefully acknowledge the hospitality  of the {\it Institut 
de Math\'ema\-ti\-ques de Luminy} at Marseille, where part of this 
article has been written. For useful comments and remarks our thanks go to Tom
Braden.

\bigskip\medskip\goodbreak
\centerline{\XIIbf 0. Preliminaries}
\medskip\nobreak\noindent
For the convenience of the reader, we recall some basic notions, 
notations and constructions to be used in the sequel.

\smallskip\noindent
{\bf 0.1} Let $V$ be a real vector space of dimension~$n$, and $A\b := 
S\b (V^*)$, the symmetric algebra on the dual vector space $V^*$, i.e., 
the algebra of real valued polynomial functions on $V$. We endow~$A\b$ 
with the even grading given by $A^2 = V^*$, a convention motivated by 
equivariant cohomology, and we let $\mm := A^{>0}$ be the homogeneous 
maximal ideal of~$A\b$. For a graded $A\b$-module $M\b$, its reduction 
$$
  \overline M\,\b \;:=\; A\b/\mm \otimes_{A\b} M\b \,,
$$
modulo $\mm$ is a graded real vector space. 
\par

For a strictly convex polyhedral cone $\sigma \subset V$, we let 
$V_\sigma \subset V$ denote its linear span. In analogy to the 
definition of $A$, we consider the graded algebra
$$
  A\b_\sigma \;:=\; S\b(V_\sigma^*) \,.
$$
We usually identify its elements with polynomial functions on the 
cone~$\sigma$.

To avoid cumbersome notation, we admit graded homomorphisms even if 
they are not of degree zero.

\smallskip\noindent
{\bf 0.2} Motivated by the coarse ``toric topology'' on a toric variety 
given by torus-invariant open sets, we consider a fan $\Delta$ 
(which need not be rational) 
in~$V$  as a finite
topological space with the  subfans as open subsets. The ``affine'' fans 
$$
  \langle \sigma \rangle \;:=\; \{\sigma\} \cup \partial\sigma
  \,\preceq\, \Delta 
  \quad\hbox{with boundary fan}\quad
  \partial\sigma \;:=\; \{\tau \in \Delta \;;\; \tau \precneqq \sigma\}
$$
form a basis of the fan topology by open sets that cannot be covered by 
smaller ones. Here $\preceq$ means that a cone is a face of another cone 
or that a set of cones is a subfan of some other fan. 

Sheaf theory on that ``fan space'' is particularly simple since a 
presheaf given on the basis already ``is'' a sheaf. In particular, for 
a sheaf~$\F$ on~$\Delta$, the equality 
$$
  \F\bigl(\left<\sigma\right>\bigr) = \F_{\sigma} 
$$
of the set of sections on the affine fan $\left<\sigma\right>$ and the 
stalk at the point~$\sigma$ holds.

Furthermore, a sheaf $\F$ on~$\Delta$ is flabby if and only if each 
restriction homomorphism $\rho^\sigma_{\partial \sigma} \: 
\F(\left<\sigma\right>) \to \F(\partial \sigma)$ is surjective.
\par

In particular, we consider (sheaves of) $\A\b$-modules, where $\A\b$ is 
the {\it structure sheaf} of~$\Delta$, i.e., the graded sheaf of 
polynomial algebras determined by $\A\b(\left<\sigma\right>) := 
A\b_\sigma$, the restriction  homomorphism $\rho^\sigma_{\tau} 
\: A\b_\sigma \to A\b_\tau$ being the restriction of functions 
on~$\sigma$ to the face $\tau \preceq \sigma$. The set of sections 
$\A\b(\Lambda)$ on a subfan $\Lambda \preceq \Delta$ constitutes the 
algebra of conewise polynomial functions on the support $|\Lambda|$ 
in a natural way.
\par

Given a homomorphism $\phi \: \F \to \F'$ of sheaves on~$\Delta$ and a 
subfan~$\Lambda$, we often write
$$
  F\b_\Lambda  \;:=\;  \F\b(\Lambda)\;,\quad
  F\b_\sigma  \;:=\; \F\b(\left<\sigma\right>)  
  \quad\hbox{and}\quad 
  \phi_{\Lambda}: F_{\Lambda} \to F'_{\Lambda}\; .
$$
Similarly, for a pair of subfans $(\Lambda, \Lambda_0)$ with 
$\Lambda_0 \preceq \Lambda$, we define
$$
  F\b_{(\Lambda, \Lambda_0)}  \;:=\;
  \ker (\rho^\Lambda_{\Lambda_{0}} \:
  F\b_\Lambda \to F\b_{\Lambda_0})\,,
$$
the submodule of sections on $\Lambda$ vanishing on $\Lambda_0$. In 
particular, for a purely $n$-dimensional subfan~$\Lambda$ (see~0.4 
below), we consider the case that $\Lambda_{0}$ is the {\sl boundary fan} 
$\partial \Lambda$, i.e., the subfan generated by those ($n{-}1$)-cones 
which are a facet of exactly one $n$-cone in $\Lambda$. The sections 
vanishing on $\partial \Lambda$ may be looked at as sections ``with 
compact support''.
\par

\smallskip\noindent
{\bf 0.3}  Let $f\:V \to W$ be a linear map inducing a map of fans 
between a fan $\Delta$ in~$V$ and a fan $\Lambda$ in~$W$, i.e., it maps 
each cone of $\Delta$ into a cone of~$\Lambda$. Let~$\A\b$ and~$\B\b$ 
denote the corresponding sheaves of conewise polynomial functions, and 
let $\F\b$ on~$\Delta$ and $\G\b$ on~$\Lambda$ be sheaves of graded 
$\A\b$- or $\B\b$-modules, respectively. For cones $\sigma \in \Delta$ 
and $\tau \in \Lambda$ with $f(\sigma) \subset \tau$, there is an 
induced homomorphism $B\b_{\tau} \to A\b_{\sigma}$ and thus, the 
structure of a $B\b_{\tau}$-module on $F\b_{\sigma}$.

\item{(a)} The {\it direct image\/} $f_*(\F\b)$ on $\Lambda$ is the 
$\B\b$-module sheaf defined 
by
$$
  f_*(\F\b)_\tau \;:=\; F\b_{f^{-1}(\tau)}
  \quad\hbox{with}\quad f^{-1}(\tau) \;:=\;
  \{ \sigma \in \Delta; f(\sigma) \subset \tau \}
  \;\preceq\; \Delta \;.
$$
The direct image of a flabby sheaf is again flabby. 

\item{(b)} The {\it inverse image\/} $f^*(\G\b)$ on $\Delta$ is the 
$\A\b$-module sheaf determined by
$$
  f^*(\G\b)_\sigma \;:=\; 
  A\b_\sigma \otimes_{B\b_\tau} G\b_\tau 
  \quad\hbox{for $\sigma \in \Delta$ and the minimal 
  $\tau \in \Lambda$ with 
  $f(\sigma) \subset \tau$.}
$$

\medskip
We are especially interested in the following maps of fans:

\item{(i)} For a subdivision $\hat\Delta$ of~$\Delta$, the mapping 
$\id_{V} \: (V, \hat\Delta) \to (V, \Delta)$ is such a morphism of 
fans. In particular, we will consider the case of an affine fan 
$\left<\sigma\right>$ given by an $n$-dimensional cone, and its stellar 
subdivision  
$$
  \hat\sigma := \partial \sigma + \lambda :=
  \{ \tau + \lambda; \tau \in \partial \sigma \} 
$$ 
with respect to a ray $\lambda := \ell \cap \sigma$, where~$\ell$ is a 
one-dimensional linear subspace passing through the interior of~$\sigma$.

\item{(ii)} For a cone $\sigma \in \Delta$, its  closure in the fan 
topology is the star 
$$
  \Delta_{\succeq \sigma} \;:=\; 
  \{\gamma \in \Delta \;;\; \gamma \succeq \sigma\} \,. 
$$
In general, it is not a fan. The collection 
$$
  \Delta/\sigma \;:=\; \{\pi(\gamma) \;;\; 
  \gamma \in \Delta_{\succeq\sigma}\}
$$ 
of the image cones with respect to the projection $\pi\: V \to W := 
V/V_{\sigma}$, however, is a fan, called the {\it transversal fan\/} 
of~$\sigma$ with respect to~$\Delta$. The induced map $\pi_{\sigma} : 
\Delta_{\succeq\sigma} \to \Delta/\sigma$ is a homeomorphism.

\item{(iii)} Applying~(ii) to the case of $\hat\sigma$ from (i), 
the projection $\pi:V \to W:=V/\ell$ maps the boundary fan 
$\partial \sigma$ homeomorphically onto the ``flattened boundary fan'' 
$$
  \Lambda_{\sigma} \;:=\; \hat\sigma /\lambda \;\cong\; \partial \sigma
$$
in~$W$. In that situation, choosing a linear form $T \in A^2$ with 
$T|_{\lambda} > 0$, we obtain isomorphisms $\ker(T) \cong W$ and thus 
$A\b \;\cong\; B\b[T]$, where we identify $B\b := S(W^*)$ with the 
subalgebra $\pi^{*}(B\b) \subset A\b$. -- Moreover, for a sheaf~$\F\b$ 
on $\left<\sigma\right>$ and  $\G\b \;:=\; 
\pi_{*}(\F\b|_{\partial\sigma})$, there is a natural isomorphism of 
$B$-modules
$$
  G_{\Lambda_{\sigma}} \;\cong\; F_{\partial \sigma}\,. 
  \leqno\hbox{(0.3.1)}
$$

\medskip\noindent
{\bf 0.4}  We use notations as $\Delta^d := \{\gamma \in \Delta \;;\; 
\dim\gamma=d\}$, $\Delta^{\le d}$, etc. The fan $\Delta$ in $V$ is 
called: 

\item{(a)} {\it  oriented\/} if for each cone $\sigma \in \Delta$, an 
orientation $\or_{\sigma}$ of $V_\sigma$ is fixed in such a way that 
orientations for full-dimensional cones coincide, 

\item{(b)} {\it  purely $n$-dimensional\/} if each maximal cone 
of~$\Delta$ lies in $\Delta^n$,  

\item{(c)} {\it irreducible\/} if it is not the union of two proper 
subfans with intersection included in $\Delta^{\le n{-}2}$, 

\item{(d)} {\it normal\/} if it is purely $n$-dimensional and for each 
cone $\sigma \in \Delta$, the transversal subfan $\Delta/\sigma$ is 
irreducible in $V/V_{\sigma}\,$, or, equivalently, if the support 
$|\Delta|$ is a normal pseudomanifold,

\item{(e)} {\it quasi-convex\/} if it is purely $n$-dimensional and the 
support $|\partial\Delta|$ of its boundary subfan is a real homology 
manifold. Note that a quasi-convex fan is normal, but not vice versa.
\bigskip

\bigskip\medskip\goodbreak
\centerline{\XIIbf 1. Pure Sheaves on a Fan}
\medskip
\noindent
We recall the definition of the class of ``pure'' sheaves on a fan space that
plays a key role in the sequel.

\medskip
\noindent
{\bf 1.1 Definition.} {\it A {\bfit pure sheaf\/} on a fan $\Delta$ 
is a {\it flabby\/} sheaf~$\F\b$ of graded $\A\b$-modules such that, for 
each cone $\sigma \in \Delta$, the $A\b_{\sigma}$-module $F\b_{\sigma} 
= \F\b(\langle \sigma \rangle)$ is finitely generated and free.}

\smallskip
We collect some useful facts about these sheaves, proved in [BBFK$_2$] 
and [BreLu$_1$]: 

Pure sheaves are built up from simple objects that correspond to 
the cones of the fan, or, equivalently, to the stalks of the structure 
sheaf. Up to a shift, such a simple sheaf is obtained from that stalk 
by a minimal extension process. 

\medskip\noindent
{\bf 1.2 Simple Pure Sheaves.\ } For each cone $\sigma \in \Delta$, 
we construct a ``minimal'' pure sheaf $\L\b =: {}_{\sigma}\L\b$ 
supported on the star $\Delta_{\succeq\sigma}$ and with stalk 
$L\b_{\sigma} = A\b_{\sigma}$ as follows: On the subfan 
$\Delta \setminus \Delta_{\succeq\sigma}$, we set $\L\b := 0$. 
By induction on the dimension, we extend it to the cones in 
$\Delta_{\succeq\sigma}$, starting with
$$
  L\b_{\sigma} \;:=\; A\b_{\sigma}\,.
$$
For a cone $\gamma \succeq \sigma$, we may assume that 
$L\b_{\partial\gamma}$ has been defined, and then set
$$
  L\b_{\gamma} \;:=\; 
  A\b_{\gamma} \otimes_{\RR} \quer{L}\b_{\partial\gamma}\,.
$$
The restriction homomorphism $\rho^{\gamma}_{\partial\gamma}$ is defined 
by the following commutative diagram
$$
  \matrix{
  L\b_{\gamma} := 
  A\b_{\gamma} \otimes_{\RR} \quer{L}\b_{\partial\gamma} & 
  \longto &
  \quer{L}\b_{\partial\gamma} \cr
  \noalign{\vskip 6pt}
  \kern-60pt\downarrow\rlap{$\rho := \id \otimes s$} &
  \swarrow & \Vert \cr
  \noalign{\vskip 6pt}
  L\b_{\partial\gamma} = 
  A\b_{\gamma} \otimes_{A\b_{\gamma}} L\b_{\partial\gamma} &
  \longto &
  \quer{L\b}_{\partial\gamma} \cr
  }
$$
where the diagonal arrow $s \: \quer L\b_{\partial\gamma} \to
L\b_{\partial\gamma}$ is an $\RR$-linear section of the reduction map 
in the bottom row. 

\medskip\noindent
{\bf 1.3 Remarks.} i) For each cone $\sigma \in \Delta$, the 
corresponding simple sheaf $\L\b := {}_{\sigma}\L\b$ is pure; it is 
characterized by the following properties: 
\item{a)} $\quer L\b_{\sigma} \cong \RR\b$,
\item{b)} for each cone $\tau \ne \sigma$, the reduced restriction 
homomorphism $\quer L\b_{\tau} \to \quer L\b_{\partial\tau}$ is an 
isomorphism.
\par

\noindent
In particular, property b) implies the vanishing of $\L\b := 
{}_{\sigma}\L\b$ outside of the star of~$\sigma$.

\smallskip
\noindent
ii) 
For the zero cone~$o$, the ``generic point'' of~$\Delta$, the 
corresponding simple sheaf
$$
  \E\b \;:=\; {}_{\Delta}\E\b \;:=\; {}_{o}\L\b
$$
is called the {\it equivariant intersection cohomology sheaf\/} (or the 
minimal extension sheaf) of~$\Delta$. For a quasi-convex fan~$\Delta$, 
we may define its (virtual) {\it intersection cohomology\/} as
$$
  IH\b (\Delta) \;:=\; \overline E\b_\Delta\; .
  \leqno(1.3.1)
$$

\noindent
iii) By extending scalars, each ``local'' sheaf ${}_{\sigma}\L\b$ is 
derived from the ``global'' sheaf ${}_{\Delta/\sigma}\E\b$ of the 
corresponding transversal fan: As in 0.3 (ii), we let $\pi_{\sigma} 
:= \pi|_{\Delta_{\succeq\sigma}} : \Delta_{\succeq\sigma} \to 
\Delta/\sigma$ denote the homeomorphism induced from the projection 
$V \to V/V_{\sigma}$. The inverse image 
$\pi_{\sigma}^{*}({}_{\Delta/\sigma}\E\b)$ is a flabby sheaf of 
graded $\A\b$-modules on the closed subset $\Delta_{\succeq\sigma}$ 
of~$\Delta$. Its trivial extension to the whole fan space~$\Delta$ then 
yields the sheaf~${}_{\sigma}\L\b$.

\medskip
The following elementary decomposition theorem has been proved in 
[BBFK$_2$, 2.4] and in [BreLu$_1$, 5.3]:

\medskip\noindent
{\bf 1.4 Decomposition Theorem.}
{\it Every pure sheaf $\F\b$ on~$\Delta$ admits a natural direct sum 
decomposition of $\A\b$-modules 
$$
  \F\b \;\;\cong\;\;
  \bigoplus_{\sigma \in \Delta}
  \bigl({}_{\sigma}\L\b \otimes_{\RR} K\b_\sigma\bigr)
$$
with $K\b_\sigma  :=  K\b_\sigma(\F\b)  :=
\ker\,(\,\quer \rho^{\sigma}_{\partial\sigma} \: \quer F\b_{\sigma} \to
\quer F\b_{\partial\sigma})$, a finite dimensional graded vector space.}

\smallskip
For a proof of the following application, we refer to [BBFK$_2$, 2.5].

\medskip
\noindent
{\bf 1.5 Example.} Let $\phi := \id_{V} : (V, \hat\Delta) \to 
(V, \Delta)$ be a refinement. Then $\phi_*(\hat \E\b)$ is a pure 
sheaf. Its decomposition is of the form 
$$
  \phi_*({\hat \E}\b) \;\;\cong\;\;
  \E\b \oplus \kern -3pt\bigoplus_{\sigma \in
  \Delta^{\ge2}} \kern -3pt
  \bigl( {}_{\sigma}\L\b \otimes_{\RR} K\b_\sigma
  \bigr) \,,
$$
where the~$K\b_\sigma$ now are (strictly) positively graded vector 
spaces and the ``correction terms'' are supported on the closed subset 
$\Delta^{\ge2}$. 

\medskip
\noindent
{\bf 1.6 Remark.} For a pure sheaf~$\F\b$ on the boundary fan 
$\partial\sigma$ of an $n$-dimensional cone and the projection mapping 
$\pi \: (V, \partial\sigma) \to (W, \Lambda_{\sigma})$ corresponding to 
a ray $\lambda$ as in 0.3,~(ii), the direct image 
$\pi_*(\F|_{\partial\sigma})$ is a pure sheaf on~$\Lambda_{\sigma}$.\qed

\bigskip
\bigskip\medskip\goodbreak
\centerline{\XIIbf 2. The Dual of a Pure Sheaf}

\medskip\nobreak\noindent
In this section, the symbol $\F\b$ always denotes a pure sheaf on 
an oriented fan~$\Delta$. Furthermore, unless otherwise stated, the 
symbol $\Hom$ is understood to mean $\Hom_{A\b}$, and $\otimes$ 
means~$\otimes_{\RR}$. Moreover, for a cone $\sigma \in \Delta$, we 
consider $\det V^*_\sigma:= \bigwedge^{\dim \sigma} V^*_\sigma$ as a 
graded vector space concentrated in degree $2 \dim \sigma$, with the 
convention $\det V_{o}^{*} = \RR$.

To $\F\b$, we associate its dual $\D\F\b$ and show the following 
properties: The dual is again a pure sheaf on~$\Delta$, and for each 
normal subfan $\Lambda$, the module of sections $(\D\F)_\Lambda$ 
is the dual of the module 
$F\b_{(\Lambda, \partial\Lambda)}$ of sections with compact supports 
of~$\F\b$. 

%
\medskip\nobreak\noindent
{\bf 2.1 Construction of the dual sheaf.\ } 
To construct the dual $\D\F\b$ of the pure sheaf $\F\b$ 
on~$\Delta$, we first define its sections over affine fans in such a 
way that duality holds by definition.

\smallskip\noindent
{\it Sections over a cone $\sigma \in \Delta$.\ } 
As $A\b_{\sigma}$-module, we define $(\D\F)_\sigma = 
\D\F\b(\langle \sigma \rangle)$ by 
$$
  (\D\F)\b_\sigma \;\;:=\;\; 
  \Hom (F\b_{(\sigma, \partial \sigma)}, A\b_\sigma) 
  \otimes \det V_\sigma^* \;. 
  \leqno(2.1.1)
$$

\smallskip\noindent
{\it Restriction homomorphisms.} The homomorphism $\rho^\sigma_{\tau}$ 
for $\sigma \succeq \tau$ is constructed in two steps: In the first 
step, we deal with  the case of a facet; in the second step, we extend 
this recursively to the general situation of a face of arbitrary 
codimension.

To that end, we need {\it transition coefficients\/} 
$\epsilon^{\sigma}_\tau = \pm 1$ for the facets $\tau$ of $\sigma$: For 
$d := \dim\sigma$, there exists a natural map $\kappa \: 
\bigwedge^{d-1}V_\sigma^* \to \bigwedge^{d-1}V_\tau^* = \det V_\tau^*$. 
We choose a  linear form~$h$ on $V_\sigma^*$ with $V_\tau= \ker (h)$ 
and $h|_\sigma \ge 0$. Every element of $\det V_{\sigma}^*$ decomposes 
in the form $h \,{\wedge}\, \eta$ with unique image $\kappa(\eta)$. We 
thus obtain a homomorphism 
$$
  \psi_h : \det V_\sigma^* 
  \;\longto\; \det V_\tau^* \;,\quad 
  h \wedge \eta \;\longmapsto\; \kappa(\eta) \;.
  \leqno (2.1.2)
$$ 
If now $\omega_\sigma \in \det V_{\sigma}^*$ and $\omega_\tau \in 
\det V_{\tau}^*$ define the orientations of $\sigma$ respectively 
$\tau$, we set 
$$
  \epsilon^\sigma_\tau := 
  \cases{+1 & if $\psi_h(\omega_{\sigma}) \in 
  \RR_{>0}\;\omega_{\tau}\,$,\cr
  \noalign{\vskip3pt}
  -1 & otherwise.\cr}
  \leqno(2.1.3)
$$

\medskip\nobreak\noindent
{\it Step 1: Restriction homomorphism for a facet 
$\tau \prec_{1} \sigma$. \ } Using again the linear form 
$h \in V_\sigma^*$, we are going to define another homomorphism 
$$
  \phi_h: \Hom ( F\b_{(\sigma, \partial \sigma)}, A\b_\sigma) 
  \;\longto\; 
  \Hom ( F\b_{(\tau, \partial \tau)}, A\b_\tau)  \leqno(2.1.4)
$$
and see that 
$$
  \phi_{\lambda h}=\lambda \phi_h 
  \quad\hbox{and}\quad
  \psi_{\lambda h}=\lambda^{-1} \psi_h
  \leqno (2.1.5)
$$
for every non-zero scalar $\lambda \in \RR$. Thus the homomorphism 
$$
  \phi_h {\otimes} \psi_h : 
  \Hom ( F\b_{(\sigma, \partial \sigma)}, A\b_\sigma) 
  \otimes \det V_\sigma^*
  \;\;\longto\;\;
  \Hom ( F\b_{(\tau, \partial \tau)}, A\b_\tau) \otimes \det V_\tau^*\,
$$
does not depend on the special choice of~$h$, and we may set
$$
  \rho^{\sigma}_\tau \;\;:=\;\; 
  \epsilon^{\sigma}_\tau \cdot \phi_h {\otimes} \psi_h \; .
$$

The map~$\phi_h$ associates to a homomorphism $f \:
F\b_{(\sigma, \partial \sigma)} \to A\b_\sigma$ the homomorphism 
$\phi_h (f) \: F\b_{(\tau, \partial \tau)} \to A\b_\tau$, which acts in 
the following way: We first extend a section $s \in 
F\b_{(\tau, \partial \tau)}$ trivially to $\partial \sigma$ and then to 
a section $\check{s} \in F\b_\sigma$; we thus have 
$h\check s \in F\b_{(\sigma, \partial \sigma)}$ and may finally set 
$$
  \phi_h (f)(s) \;:=\; f(h \check s)|_\tau \in A_\tau \,.
  \leqno(2.1.6)
$$
In order to see that this definition is independent of the particular 
choice of $\check s$, we  present an alternative description, following 
the argument on p.36 of [BBFK]:  We use three exact sequences, starting 
with 
$$
  0 \wideto F\b_{(\sigma,\partial \sigma)} \longto
  F\b_\sigma \longto  F\b_{\partial \sigma} \longto 0 \;.
  \leqno(2.1.7)
$$
The second one is composed of the multiplication with $h$ and the 
projection onto the cokernel:
$$
  0 \longto A\b_\sigma \;{\buildrel \mu_h \over \longto}\;
  A\b_\sigma  \longto A\b_\tau \longto 0\;.
  \leqno(2.1.8)
$$
Eventually the subfan $\partial_\tau \sigma:= \partial\sigma \setminus 
\{\tau \}$ of $\partial\sigma$  gives rise to the exact sequence 
$$
  0 \wideto F\b_{(\tau, \partial \tau)}
  \wideto F\b_{\partial \sigma}
  \wideto F\b_{\partial_\tau  \sigma}
  \wideto 0\;.\leqno(2.1.9)
$$
The associated $\Hom$-sequences provide a diagram
$$
\def\longto{\kern-6pt\longrightarrow\kern-6pt}
\diagram{
    &   &   &   & \Ext(F_{\partial_\tau\sigma}\b,
A\b_\sigma)
                    &   &   \cr 
\noalign{\leftline{(2.1.10)}\vskip-20pt}             
    &   &   &   & \mapdown{}
                    &   &   \cr 
  \Hom(F\b_\sigma, A\b_\sigma) & \longto &
          \Hom(F\b_{(\sigma,\partial\sigma)}, A\b_\sigma) &
              \buildrel \alpha\over\longto &
                  \Ext(F\b_{\partial\sigma}, A\b_\sigma)
                    &   &   \cr 
    &   &   &   & \mapdown{\beta}
                    &   &   \cr 
  \Hom( F_{(\tau,\partial\tau)}\b, A\b_\sigma ) & \longto &
          \Hom( F_{(\tau,\partial\tau)}\b, A_\tau\b) &
              \buildrel \gamma\over\longto &
                  \Ext(F_{(\tau,\partial\tau)}\b, A\b_\sigma)
                    & \longto &
\Ext(F_{(\tau,\partial\tau)}\b, A\b_\sigma)
\cr 
}
$$
with $\Ext=\Ext^1_{A\b}$. We show that~$\gamma$ is an isomorphism; we 
then have
$$
  \phi_h \;:=\; \gamma^{-1}\circ\beta\circ\alpha \,.
  \leqno (2.1.11)
$$
In fact, the rightmost arrow in the bottom row is the zero homomorphism, 
since it is induced by multiplication with $h$, which annihilates 
$F\b_{(\tau, \partial \tau)}$. On the other hand, the $A\b_\tau$-module 
$F\b_{(\tau, \partial \tau)}$ is a torsion module over $A\b_\sigma$, so 
that $\Hom(F_{(\tau,\partial\tau)}\b, A\b_\sigma)$ vanishes.
\par

\medskip\nobreak\noindent
{\it Step 2: Restriction homomorphism for faces of higher codimension.\ } 
For a face $\tau \prec \sigma$ of codimension~$r \ge 2$, we choose a 
``flag'' 
$$
  \tau =: \tau_0 \;\prec_1\; \tau_1 \;\prec_1\; 
  \ldots \;\prec_1\; \tau_r := \sigma
$$ 
of relative facets joining $\tau$ and $\sigma$. Defining the restriction 
homomorphism $\rho^\sigma_{\tau}$ as the composition of the 
$\rho^{\tau_{i+1}}_{\tau_{i}}$, we have to show that the result does not 
depend on the particular choice of the flag. This is easy to see in the 
case $r=2$: For two flags $\gamma \prec_1 \tau \prec_1 \sigma$ and 
$\gamma \prec_1 \tau' \prec_1 \sigma$ and $h, h' \in V^*_\sigma$ as 
above, we set $g := h|_{V_{\tau'}}$ and $g' := h'|_{V_{\tau}}$, and then 
find
$$
  \phi_{g'} \circ \phi_{h} \;=\;
  \phi_{g} \circ \phi_{h'} \quad \hbox{and} \quad 
  \psi_{g'} \circ \psi_{h} \;=\;
  - \kern1.5pt \psi_{g} \circ \psi_{h'}\ ,
$$
whence
$$
  \rho^{\tau}_\gamma \circ \rho^{\sigma}_\tau
  \;=\; \rho^{\tau'}_\gamma \circ \rho^{\sigma}_{\tau'}\;.
  \leqno(2.1.12)
$$

Thus, for general~$r$, it suffices to verify that every two such flags 
can be transformed into each other in such a way that in each step, 
only one intermediate cone is replaced by another one. We proceed by 
induction on the codimension~$r$. 

To prove that claim, we may assume $\tau = o$ (otherwise, we replace 
$\Delta$ with $\Delta/\tau$) and $\Delta = \langle \sigma \rangle$. We 
want to compare the given flag with a second one, say 
$o \prec_1 \tilde \tau_{1} \prec_1 \tilde \tau_2 \prec_1 \ldots \prec_1
\tilde \tau_r=\sigma$. There is a chain of rays $\rho_{1} := \tau_{1}, 
\ldots, \rho_s := \tilde \tau_{1}$ such that the two-dimensional cones 
$\rho_i{+}\rho_{i+1}$ belong to~$\Delta$. We now proceed by a second 
induction on~$s$. For $s=1$, we may pass to the fan $\Delta/\tau_{1}$ 
and use the first induction hypothesis for $r{-}1$. For the induction 
step, it evidently suffices to consider the case $s=2$. Choosing any 
auxiliary flag of the form $o \prec_{1} \tau_{1} \prec_{1} 
\tau_{1}{+}\tilde\tau_{1} \prec_{1} \ldots \prec_{1} \sigma$, the case 
$s=1$ yields its equivalence with the start flag. On the other hand, by 
(2.1.12), the auxiliary flag is equivalent to the one obtained by 
interchanging $\tau_{1}$ and $\tilde\tau_{1}$, and this in turn is 
equivalent to the ``twiddled'' flag.\qed


We now show that the formula (2.1.1) for cones extends even to normal 
subfans. To that end, we need the following preparatory results:

\bigskip
\noindent
{\bf 2.2 Lemma.\ } 
{\it{\rm(i)}  For an arbitrary fan~$\Delta$, there is a natural 
isomorphism
$$
  \Theta:
  \bigoplus_{\sigma \in \Delta^n} (\D\F)_{\sigma}
  \;\buildrel \cong \over \longto\;
  \Hom (F\b_{(\Delta,\Delta^{\le n-1}
  )},  A\b) \otimes \det V^*\; .
  \leqno(2.2.1)
$$
{\rm(ii)} If $\Delta$ is purely $n$-dimensional, the $A$-modules
$$
  F_{(\Delta, \Delta^{\le n-1})} \subset F_{(\Delta, \partial \Delta)} 
  \subset F_\Delta\, 
$$
are torsion-free and of the same rank. 
As a consequence, the restriction homomorphisms 
$$
  \Hom(F_{\Delta}, A) \longto 
  \Hom(F_{(\Delta, \partial \Delta)}, A) \longto 
  \Hom(F_{(\Delta, \Delta^{\le n-1})}, A)
$$ 
for the dual modules are injective.

{\rm(iii)} In the setup of cellular  (``{\v C}ech'') 
cochains and cocycles as in section~3 of\/ {\rm[BBFK${}_2$]}, for 
an arbitrary sheaf~$\G$ on~$\Delta$, there is an isomorphism 
$$
  Z^0(\Delta; \G) = G_{(\Delta, \partial \Delta)}\, ,
  \leqno{\rm(2.2.2)}
$$
and for a normal fan, we also have an isomorphism}
$$
  Z^0(\Delta, \partial\Delta; \G) = G_{\Delta}\,. 
  \leqno{\rm(2.2.3)}
$$

\medskip
\noindent
{\bf Proof:} (i) For each $n$-dimensional cone $\sigma$, the equality 
$\det V_{\sigma}^* = \det V^*$ holds. Hence, the isomorphism $\Theta$ 
is immediately obtained from the defining equality (2.1.1) by applying 
the additive functor $\Hom(\_\,, A)\otimes V^*$ to the obvious direct 
sum decomposition
$$
  F_{(\Delta, \Delta^{\le n-1})}
  \cong \bigoplus_{\sigma  \in \Delta^n} 
  F_{(\sigma, \partial \sigma)}\ . 
$$
(ii) For the special case $\F\b = \E\b$, this has been proved in 
[BBFK$_{2}$, 6.1, i)]. The proof clearly carries over to arbitrary 
pure sheaves.

\noindent (iii) We recall that the submodule 
$$
  Z^0 (\Delta, \partial \Delta; \G) \subset
  C^0 (\Delta, \partial \Delta; \G) = 
  \bigoplus_{\sigma \in \Delta^n} G_\sigma
$$
of degree zero cocycles $g = (g_\sigma)$
relative to $\partial \Delta$ consists of those 
cochains that satisfy 
$g_\sigma|_\tau = g_{\sigma'}|_\tau$ 
whenever $\tau \in \Delta^{n-1}$ is a common 
facet of two $n$-cones $\sigma$ and $\sigma'\,$, 
whereas for the submodule 
$$
  Z^0 (\Delta; \G) \subset
  Z^0 (\Delta, \partial \Delta; \G)
$$
of absolute cocycles, we require in addition
that the restriction of $g_\sigma$ to each ``outer'' 
facet $\tau \in \partial \Delta^{n-1}$ vanishes.

For (2.2.2) and (2.2.3)  
we note that the right hand side is always contained 
in the left hand side. In order to see the reverse inclusion, 
we have to show that $g_\sigma|_\gamma = g_{\sigma'}|_\gamma$ 
holds whenever $\gamma \in \Delta$ is a common face of two 
$n$-cones $\sigma$ and $\sigma'$. \TOT
For an ``outer'' facet 
$\gamma \in \partial \Delta$ and an absolute cocycle~$g$, 
that is clear, while for an ``inner'' facet, we may use that 
$\Delta/\gamma$ is irreducible: In the case 
$\gamma \not\in \partial \Delta$, the fan $\Delta/\gamma$ is 
complete and thus irreducible, while the case 
$\gamma \in \partial \Delta$ only has to be considered for a relative
cocycle and then it follows from the fact that $\Delta$ is assumed 
to be normal. Finally we use that in an irreducible fan, two cones 
of maximal dimension may be connected by a chain of cones such 
that successive cones have a facet in common.
\LEB
Since $\Delta$ is normal, the cones $\sigma$ and $\sigma'$ can be joined 
by a chain of $n$-cones intersecting successively in common facets containing
$\gamma$. It thus suffices to consider the case that $\gamma$ is a facet,
where the statement is obvious.\qed

\bigskip
\noindent
{\bf 2.3 Theorem.\ }
{\it For a normal oriented fan~$\Delta$, the natural isomorphism $\Theta$ 
of\/ {\rm(2.2.1)} induces isomorphisms}
$$
  (\D\F)_\Delta
  \;\buildrel \cong \over \longto\;
  \Hom (F\b_{(\Delta, \partial \Delta)},  A\b) \otimes \det V^*
  \leqno(2.3.1)
$$
{\it and}
$$
  (\D\F)_{(\Delta, \partial \Delta)}
  \;\buildrel \cong \over \longto\;
  \Hom (F\b_{\Delta},  A\b) \otimes \det V^*\; .
  \leqno(2.3.2)
$$

\medskip
\noindent
{\bf Proof:} \TOT{For a sheaf $\G$ on $\Delta$ we denote 
$$
  Z^0 (\Delta, \partial \Delta; \G) \subset
  \bigoplus_{\sigma \in \Delta^n} G_\sigma
$$
the submodule of cocycles $g= (g_\sigma)$ relative 
$\partial \Delta$, i.e. we require that $g_\sigma|_\tau = 
g_{\sigma'}|_\tau$ holds whenever $\tau \in \Delta^{n-1}$ is a common 
facet of two $n$-cones $\sigma$ and $\sigma'\,$, and
$$
  Z^0 (\Delta; \G) \subset
  Z^0 (\Delta, \partial \Delta; \G)
$$
the submodule of absolute cocycles, i.e. we require in addition that 
the restriction of $g_\sigma$ to each ``outer'' facet 
$\tau \in (\partial \Delta)^{n-1}$ vanishes. In fact, 
$$
  Z^0(\Delta; \G) = G_{(\Delta, \partial \Delta)}\ , 
$$
while for a normal fan we  even have
$$ 
  Z^0(\Delta, \partial \Delta; \G) =
  G_\Delta\; . 
$$
In both cases the right hand side is contained in the left hand side; 
in order to see the reverse inclusion we have to show that 
$g_\sigma|_\tau = g_{\sigma'}|_\tau$ holds whenever $\tau \in \Delta$ 
is a common facet of two $n$-cones $\sigma$ and $\sigma'$. For $\tau \in
\partial \Delta$ and an absolute cocycle $g$, that is clear, while 
otherwise we may use that $\Delta/\tau$ is irreducible: In the case 
$\tau \not\in \partial \Delta$, the fan $\Delta/\tau$ is complete and 
thus irreducible, while the case $\tau \in \partial \Delta$ only has to 
be considered for a relative cocycle, and then it follows from the fact 
that $\Delta$ is assumed to be normal. Finally use that in an 
irreducible fan two cones of maximal dimension may be connected by a 
chain of cones such that successive cones have a facet in common.}
\LEB
Using the formalism of \v{C}ech cochains as in (iii) of the Lemma, we 
restate the
assertion as follows: For any 
$0$-cochain $\psi = (\psi_\sigma) \in \bigoplus_{\sigma \in \Delta^n} 
\D\F_\sigma$ and its image $\Theta (\psi)$ in 
$\Hom (F\b_{(\Delta,\Delta^{\le n-1})},  A\b) \otimes \det V^*$, the 
equivalences 
$$
  \psi \in Z^0(\Delta, \partial \Delta; \D\F) \iff 
  \Theta (\psi) \in 
  \Hom (F_{(\Delta, \partial \Delta)}, A) \otimes \det V^* 
  \leqno (2.3.3)
$$
and
$$
  \psi \in Z^0(\Delta; \D\F) \iff
  \chi =\Theta (\psi) \in \Hom (F_{\Delta}, A) \otimes \det V^* 
  \leqno (2.3.4) 
$$ 
hold. 

To prove these equivalences, we choose an auxiliary function 
$h = \prod_{i=1}^r h_i \in A^{2r}$ as the lowest degree product of 
linear forms~$h_{i}$ that vanishes on $\bigcup_{\tau \in \Delta^{n-1}} 
V_{\tau}$; so each $V_\tau$ is the kernel~$V_{i}$ of some~$h_{i}$. 
After fixing a positive volume form on~$V$ and thus, an isomorphism 
$\RR \buildrel \cong \over \longto \det V^{*}$, the homomorphisms 
$\psi_{h_{i}}$of (2.1.2) provide isomorphisms $\RR \cong \det V^{*}
\cong \det V_{i}^{*}$. We may thus drop the determinant factors on the 
right hand side, and for each cone $\gamma \in\Delta^{\ge n-1}$, we may 
replace $(\D\F)_\gamma$ with $\Hom (F\b_{(\gamma, \partial \gamma)},
A_{\gamma}\b)$ and the restriction maps with $\pm \phi_{h_i}$. Using 
the obvious inclusions 
$$
  hF_{(\Delta, \partial \Delta)} \subset 
  h F_\Delta \subset F_{(\Delta, \Delta^{\le n-1})}
$$
of torsion-free $A$-modules, the right hand sides of (2.3.3) and (2.3.4) 
are equivalent to the inclusions 
$$
  \chi (hF_{(\Delta, \partial \Delta)}) 
  \subset hA \quad\hbox{ and }\quad\chi (hF_\Delta) \subset hA ,
  \quad\hbox{where}\quad \chi:=\Theta(\phi)\;.
$$ 

\smallskip\noindent
``$\Longrightarrow$": In order to prove these implications, it suffices 
to show that for a pertinent $0$-cocycle~$\psi$, the divisibility 
relation $h_i \,|\, \chi (hf)$ holds for each index~$i$ and for an 
arbitrary section~$f$ in $F_{(\Delta, \partial \Delta)}$ or $F_\Delta$, 
respectively. 

With $f_\sigma := f|_\sigma \in F_\sigma$ for an $n$-cone~$\sigma$, we 
write 
$$
  \chi (hf) = 
  \sum_{\sigma \in \Delta^n} \psi_\sigma (hf_\sigma) \in A\b\,. 
$$
For each index $i=1, \ldots, r$, we introduce the monomial $g_{i} := 
h/h_{i} \in A^{2r-2}$. For the implication in (2.3.3), we consider a 
relative $0$-cocycle $\psi = (\psi_{\sigma}) \in Z^0 (\Delta, 
\partial \Delta; \D\F)$ and a section $f \in 
F_{(\Delta, \partial \Delta)}$ ``with compact support''. If 
$\sigma \,{\cap}\, V_i$ is not a facet or belongs to $\partial \Delta$, 
then $g_i f_\sigma$ lies in $F\b_{(\sigma, \partial \sigma)}$ and thus 
$\psi_\sigma (hf_\sigma) = h_i \psi_\sigma(g_if_\sigma) \in h_i A\b$ 
holds.  Otherwise, there is precisely one $n$-cone $\sigma' \ne \sigma$ 
such that~$\tau := \sigma \,{\cap}\, V_i$ is a common facet of both, 
$\sigma$ and $\sigma'$. We now verify that $h_i$ divides the sum 
$\psi_\sigma (hf_{\sigma}) + \psi_{\sigma'} (hf_{\sigma'})$ or, 
equivalently, that
$$
  \psi_{\sigma}(hf_{\sigma})\,|_{V_{i}} \;=\;
  - \psi_{\sigma'}(hf_{\sigma'})\,|_{V_{i}} 
  \leqno(2.3.5)
$$
holds. Using the extension $g_{i}f_{\sigma} \in 
F_{(\sigma, \partial_\tau \sigma)}$ of $(g_{i}f)|_{\tau} \in 
F_{(\tau,\partial \tau)}$ in formula  (2.1.6), we obtain
$$
  \psi_\sigma (hf_{\sigma})|_{V_i} \;\;=\;\;
  \bigl(\phi_{h_i} (\psi_\sigma)\bigr)
  \bigl((g_if_{\sigma})|_{\tau}\bigr)\; .
$$
By the relative cocycle condition, $\psi_\sigma$ and $\psi_{\sigma'}$ 
restrict to the same section in $(\D\F)_\tau$. According to the choice 
of the transition coefficients $\epsilon^\sigma_{\tau}$ in the 
definition of the restriction homomorphism $\rho^\sigma_{\tau}$ in 
(2.1.8), that yields
$$
  \phi_{h_i} (\psi_\sigma) \;=\; 
  -\phi_{h_i}(\psi_{\sigma'})\ ,
$$
which implies our claim. If $\psi \in Z^0(\Delta; \D\F)$ is an absolute 
cocycle and $f \in F_\Delta$, the argument is as above, only in the case 
that $\tau:=\sigma \cap V_i$ is an ``outer'' facet of $\sigma$ (i.e., 
contained in $\partial \Delta$), one has to use the fact that 
$\psi_\sigma|_\tau = 0$.
\par

\smallskip\noindent
"$\Longleftarrow$": For this implication, we assume that $\chi = 
\Theta(\psi) : F_{(\Delta, \Delta^{\le n-1})} \to A$ is a homomorphism 
which can be extended to the larger modules 
$F_{(\Delta, \partial \Delta)}$ or $F_\Delta$, respectively. We have to 
show the pertinent cocycle condition for $\psi = (\psi_{\sigma})$, 
namely, the equality $\psi_\sigma|_\tau = \psi_{\sigma'}|_\tau$ whenever 
$\tau$ is a common ``inner'' facet of two $n$-cones $\sigma,\sigma' \in 
\Delta$, and in the second (``absolute'') case, the vanishing 
$\psi_\sigma|_\tau=0$ if $\tau$ is an ``outer'' facet of $\sigma$.

Let $i$ be the index with $\ker (h_i)= V_\tau$. We fix an arbitrary 
section $f_0 \in F\b_{(\tau, \partial\tau)}$ and, as for (2.1.6),  
extend it to sections $f \in F\b_\sigma$, $f' \in F\b_{\sigma'}$ 
vanishing on all the remaining facets of $\sigma$ and of $\sigma'$, 
respectively. Patching them together and extending by 0 yields a 
section $f_1 \in F\b_{(\Delta, \partial \Delta)}$. Then the equation 
$$
  h_i \chi (f_1) = \chi (h_i f_{1})
  = \chi (h_i f + h_i f') =
  \psi_\sigma (h_i f) + \psi_{\sigma'}(h_i f')\;,
$$
after restriction to $\tau$, yields
$$
  0 = \bigl(h_i \chi(f_1) \bigr)\big|_{\tau} =  
  \phi_{h_i} (\psi_\sigma) (f_0) + 
  \phi_{h_i} (\psi_{\sigma'})(f_0) = 
  (\psi_\sigma|_{\tau}  - 
  \psi_{\sigma'}|_{\tau})(f_{0})\;.
$$
Finally, we leave it to the reader to consider the remaining case where 
$\chi \in \Hom (F_\Delta, A)$ and $\tau \in \partial \Delta$.\qed

\bigskip
\noindent
{\bf 2.4 Theorem:} {\it The dual sheaf $\D\F$ of a pure sheaf $\F\b$ 
is again pure.}

\bigskip
\noindent
{\bf Proof}: As in Corollary 4.12 in [BBFK], the $A\b_\sigma$-module 
$F\b_{(\sigma, \partial \sigma)}$ is free and thus also its dual 
$(\D\F)_\sigma$; hence, we only have  to prove that, for each cone 
$\sigma \in \Delta$, the restriction {homomorphism}
$$
  \rho^\sigma_{\partial\sigma} :
  (\D\F)_\sigma \longto (\D\F)_{\partial \sigma}
$$
is surjective.

To that end, we first interpret $(\D\F)_{\partial \sigma}$. We may 
assume $\dim \sigma = n$ and use the setup of 0.3~(iii). For $\G\b := 
\pi_* (\F\b|_{\partial \sigma})$, as in (0.3.a), there is a natural 
isomorphism
$$
  (\D\G)_{\Lambda_\sigma} \;\cong\; 
  (\D\F)_{\partial \sigma}
  \leqno(2.4.1)
$$
of $B\b$-modules, while for the complete fan~$\Lambda_\sigma$ {in~$W$}, 
Theorem~2.3 yields
$$
  (\D\G)_{\Lambda_\sigma} \;\cong\;
  \Hom_{B\b} (G\b_{\Lambda_\sigma}, B\b) \otimes \det W^*\; .
$$
Using the isomorphism (0.3.1), we thus obtain a chain of isomorphisms 
$$
  (\D\F)_{\partial\sigma} \cong (\D\G)_{\Lambda_\sigma}
  \cong
  \Hom_B (G_{\Lambda_\sigma}, B\b) \widecong
  \Hom_B (F_{\partial\sigma}, B\b)\;. \leqno(2.4.2)
$$
Eventually, using these isomorphisms, a section $\beta \in 
(\D\F)_{\partial \sigma}$ may be interpreted as an element of 
$\Hom_B (F_{\partial\sigma}, B\b)$. 

To proceed with the proof, we introduce the sheaf $\H\b := \pi^*(\G\b)$ 
on~$\hat\sigma$. There are isomorphisms 
$$
  H\b_{\hat \sigma}
  \widecong A\b \otimes_{B\b} G\b_{\Lambda_\sigma}
  \widecong A\b \otimes_{B\b} F\b_{\partial \sigma}
  \quad {\rm and}  \quad
  H\b_{\partial\hat \sigma} \widecong F\b_{\partial \sigma}
  \leqno(2.4.3)
$$
and a ``Thom isomorphism''
$$
 \mu_g: H\b_{\hat \sigma} \buildrel \cong \over 
 \longto g H\b_{\hat \sigma} = 
 H\b_{(\hat \sigma, \partial \hat \sigma)}\ , \
 h \mapsto gh
 \leqno(2.4.4)
$$
with a conewise linear function $g  \in 
A^2_{(\hat \sigma, \partial \hat \sigma)}$, unique up  to a non-zero 
scalar multiple, that is constructed conewise as follows: We fix a 
nontrivial linear form $f \in A^2_\lambda$. For a facet 
$\tau \prec_{1} \sigma$, let $g_\tau \in A^2$ be the unique linear form 
with $\ker (g_\tau) = V_\tau$ and $g_\tau |_\lambda = f$. Then we set 
$g|_{\tau + \lambda} := g_\tau$. 

For each facet~$\tau$ of $\sigma$, the function $g_{\tau}$ induces an 
isomorphism
$$
  \det V^* \cong \det V^*_\tau\;,\;
  g_\tau \wedge \eta \mapsto \eta|_{V_\tau} \;.
$$
Then the composed isomorphism $\det V^* \cong \det V^*_\tau \cong 
\det W^*$ is independent of $\tau$. We thus may drop the determinant 
factors.

We want to show that an inverse 
image $\alpha \in (\D\F)_{\sigma} = 
\Hom(F\b_{(\sigma, \partial \sigma)}, A\b)$ 
of~$\beta \in (\D\F)_{\partial \sigma}$ with respect to 
$\rho^\sigma_ {\partial \sigma}$ is given by the composition
$$
  F\b_{(\sigma,\partial \sigma)}
  \buildrel i \over \longto
  H\b_{(\hat \sigma, \partial \hat \sigma)}
  \buildrel \mu_{g^{-1}} \over
  \longto H\b_{\hat \sigma} \cong 
  A\b \otimes_{B\b} F\b_{\partial \sigma}
  \buildrel {\rm id}_{A} \otimes \beta \over \longto 
  A\b \otimes_{B\b} B\b = A\b
$$
where $\mu_{g^{-1}}$ is the isomorphism ``division by~$g$'' 
corresponding to (2.4.4), and the homomorphism $i$ is constructed as 
follows: Since $F\b_{\sigma}$ is a free $A\b$-module and the restriction 
homomorphism $H\b_{\hat \sigma} \to H\b_{\partial \hat\sigma}$ is 
surjective, cf.\ (2.4.3), the operator $\rho^\sigma_{\partial\sigma}$ 
for the sheaf $\F\b$ admits a factorization of the form
$$
  F\b_{\sigma}
  \buildrel j \over \longto
  H\b_{\hat \sigma} \onto
  H\b_{\partial \hat \sigma} \cong
  F\b_{\partial \sigma} \,.
$$
Since $j(F_{(\sigma, \partial \sigma)})$ clearly is included in 
$H_{(\hat \sigma, \partial \hat \sigma)}$, we may choose $i :=  
j|_{F_{(\sigma,\partial \sigma)}}$.

To prove the equality $\alpha|_{\partial\sigma} = \beta$, it still 
remains to show that $\alpha|_\tau = \beta|_\tau$ for all facets 
$\tau \prec_1 \sigma$. Here we identify the naturally isomorphic 
algebras $B\b$ and $A\b_\tau$.

We fix an arbitrary section $s \in F\b_{(\tau, \partial \tau)} \subset 
F\b_{\partial \sigma}$, where the inclusion is given by trivial 
extension. Using the isomorphisms (2.4.3) and (2.4.4), any further 
extension $\check{s}$ of~$s$ to a section of~$\F$ on the whole cone 
$\sigma$, looked at as section in $H_{\hat\sigma} \supset F_\sigma$, 
can be written in the form 
$$
  \check{s} = 1 \otimes s + g d 
  \in H\b_{\hat \sigma} \cong 
  A\b \otimes_{B\b} F\b_{\partial\sigma}
$$ 
with some correction term $d \in H\b_{\hat \sigma}$. Recalling the 
formula (2.1.6) in the definition of the homomorphism 
$\rho^\sigma_{\tau}$ for $\D\F$, we have to show that the restriction 
of the polynomial function $\alpha(g_{\tau} \cdot \check{s}) \in 
A\b_{\sigma}$ to~$\tau$ coincides with $\beta(s)$. To that end, we note 
that $g_\tau \cdot (1 \otimes s) = g \cdot (1 \otimes s)$ holds, since 
the support of $1 \otimes s \in H\b_{\hat \sigma}$ is contained in 
$\tau + \lambda$. So we eventually have the equality 
$$
  \alpha (g_\tau \cdot \check{s})|_{\tau} = 
  (\id_A \otimes \beta) (1 \otimes s + g_\tau d)|_{\tau} = 
  (\id_A \otimes \beta) (1 \otimes s)|_{\tau}\;,
$$ 
and thus $\rho^\sigma_{\tau}(\alpha)$ maps~$s$ to $\beta (s)$.\qed


\bigskip
In order to see that the dual sheaf $\D\F$ of a simple pure sheaf 
again is simple, we need biduality:

\bigskip
\noindent
{\bf 2.5 Biduality Theorem.} [BreLu${}_1$, 6.23] {\it Every pure sheaf on an
oriented  fan is reflexive: For such a sheaf~$\F\b$, there exists a natural 
isomorphism}
$$
  \F\b \buildrel \cong \over \longto \D(\D\F)\ .
$$

\noindent
{\bf Proof\/}: Over a cone $\sigma \in \Delta$, the biduality 
isomorphism $F_\sigma \to \D\D\F_\sigma$ is obtained using these 
isomorphisms:
$$
  \eqalign {(\D\D\F)_{\sigma}  = \;\; &
  \Hom\bigl((\D\F)_{(\sigma,\partial \sigma)},
  A\b_{\sigma}\bigr) \otimes \det V_\sigma^*
  \cr
  \;\;\cong\;\; &
  \Hom\bigl(\Hom(F\b_{\sigma}, A\b_{\sigma})
  \otimes \det V_\sigma^*,
  A\b_{\sigma}\bigr)\otimes \det V_\sigma^*\; \cr
  \;\cong\;\; &
  \Hom\bigl(\Hom(F\b_{\sigma}, A\b_{\sigma})
  \otimes \det V_\sigma^*,
  A\b_{\sigma} \otimes \det V_\sigma^*\bigr) \cr
  \;\cong\;\; &
  \Hom\bigl(\Hom(F\b_{\sigma}, A\b_{\sigma}),
  A\b_{\sigma} \bigr) \cr}\ \ , 
  \leqno(2.5.1)
$$
where the first isomorphism follows from Theorem 2.3 with the fan 
$\Delta := \langle \sigma \rangle$ in the vector space $V_\sigma$. The 
free $A\b_{\sigma}$-module $F\b_{\sigma}$ is reflexive, so it can be 
naturally identified with the fourth module in (2.5.1). Since this 
conewise construction is natural, it carries over to the sheaves.
\qed

\bigskip
\noindent
{\bf 2.6 Corollary.} [BreLu${}_1$, 6.26] {\it For each cone $\sigma \in \Delta$,
the simple  pure sheaf ${}_{\sigma}\L\b$ satisfies
$$
  \D ({}_{\sigma}\L\b) \;\cong\; 
  {}_{\sigma}\L\b \otimes \det\,V_\sigma^*\;.
$$
In particular, the equivariant intersection cohomology sheaf $\E\b$ is 
self-dual with an isomorphism
$$
  \theta : \E\b \;\buildrel\cong\over\longto\; \D\E 
  \leqno{(2.6.1)}
$$
of degree zero. }

\bigskip
\noindent
{\bf Proof}: Clearly, by biduality, $\D\F = 0$ implies $\F=0$. On the 
other hand, the duality functor respects a direct sum decomposition of 
pure sheaves. Since the bidual $\D\bigl(\D({}_\sigma\L)\bigr) \cong
{}_\sigma\L$ is simple, the Decomposition Theorem 1.4 implies that the 
dual $\D({}_\sigma\L)$ must be a simple sheaf. For a pure sheaf~$\F$ 
and a cone $\sigma \in \Delta$, the $A_\sigma$-module 
$F_{\partial \sigma}$ is a torsion module, whence $F_\sigma=0$ if and 
only if $F_{(\sigma, \partial \sigma)}=0$. Hence a pure sheaf and 
its dual have the same support, so $\D{}(_\sigma\L)$ and ${}_\sigma\L$ 
agree up to a shift. To determine it explicitly, we use the equality
${}_\sigma L_{(\sigma, \partial \sigma)} = {}_\sigma L_\sigma=A_\sigma$, 
which yields $D ({}_{\sigma}\L)_\sigma \,\cong\, 
\Hom (A\b_\sigma, A_\sigma) \otimes \det\,V_\sigma^* \cong 
A_\sigma \otimes \det\,V_\sigma^*$.
\qed
\bigskip

\bigskip\medskip\goodbreak
\centerline{\XIIbf 3. The Intersection Product}
\medskip\nobreak\noindent
In order to make precise the naturality of the intersection product we 
need this notion:

\medskip
\noindent
{\bf 3.1 Definition.} {\it A  {\bfit duality correlation} on $\Delta$ 
is a sheaf homomorphism
$$
  \phi : \E \longto \D\E
$$
of degree $0$ from the equivariant intersection cohomology sheaf to its 
dual extending the natural identification}  
$$
  \E_o = \RR 
  \;\buildrel  1 \mapsto 1^* \over
  \longto\;
  \RR^* = \D\E_o\ .
$$

\medskip
After multiplication with an appropriate scalar factor if necessary, 
any isomorphism $ \E \to\; \D\E$ is such a duality correlation. 
We aim at the following result:

\medskip
\noindent
{\bf 3.2 Theorem.} {\it On every fan $\Delta$, there is a unique duality 
correlation. It defines a self duality $\E \cong \D\E$ for the 
equivariant intersection cohomology sheaf $\E$.}

\medskip
Existence has already been shown in 2.6. Before proving uniqueness, we 
first use the correlation to introduce an intersection product.

\bigskip
\noindent
{\bf 3.3 Remark and Definition.} Let $\Delta$ be a normal 
$n$-dimensional oriented fan. If we fix a positive volume form 
$\omega \in \det V^*$, then every duality correlation $\phi$ gives rise 
to an {\it intersection product\/} on~$\Delta$, i.e., a pairing
$$
  E\b_\Delta \times 
  E\b_{(\Delta, \partial \Delta)}
  \;\;\longto\;\; 
  A\b[-2n] \leqno{\rm(PD)}
$$
as follows:  The isomorphism $\Theta$ of (2.3.1) yields an isomorphism 
$$
  (\D\E)_\Delta 
  \;\buildrel {\Theta_V}
\over \longto \; 
  \Hom
  \bigl(E\b_{(\Delta, \partial \Delta)}, A\b) 
  \otimes \det\, V^*
  \;\buildrel \cong \over
  \longto\;  
  \Hom (E\b_{(\Delta, \partial \Delta)}, A\b[-2n]\bigr)
  \;;\leqno{\hbox{(D$\omega$)}}
$$
its composition with the duality correlation $\phi_{\Delta}$ on the 
level of global sections provides a homomorphism
$$
  \chi_{\Delta} := \chi_{\Delta}^\omega : 
  E\b_\Delta
  \;\; \longto \;\;
  \Hom \bigl(E\b_{(\Delta, \partial \Delta)}, 
  A\b[-2n]\bigr)\,,\leqno(3.3.1)
$$
which is equivalent to (PD). 

\medskip
\noindent
{\bf 3.4 Theorem.} [BBFK${}_2$, 6.3] and [BreLu${}_1$, 6.28] 
{\it Let the oriented fan $\Delta$ be normal, and fix a positive volume 
form $\omega \in \det V^*$. If a duality correlation $\phi : \E \to 
\D\E$ is an isomorphism, then the induced pairing
$$
  E\b_\Delta \times E\b_{(\Delta, \partial \Delta)}
  \;\;\longto\;\; A\b[-2n]\ \leqno{\rm(PD)}
$$
is a dual pairing of reflexive $A$-modules. If $\Delta$ is even 
quasi-convex, then the $A$-modules $E_\Delta$ and $E\b_{(\Delta, 
\partial \Delta)}$ are free, and thus the associated reduced pairing 
$$
  \quer E\b_\Delta \times \quer E\b_{(\Delta, \partial \Delta)}
  \;\;\longto\;\; \quer A\b[-2n]
  \cong \RR[-2n]\ .\leqno (\quer {\rm PD})
$$
is a dual pairing of graded real vector spaces.}

\medskip
\noindent
{\bf Proof\/}: Compose the isomorphisms $\phi_\Delta$ and 
$\phi_{(\Delta, \partial \Delta)}$ with the isomorphisms (2.3.1) and 
(2.3.2):
$$
\displaylines{
  E_\Delta \;\buildrel \cong \over
  \longto\; \D\E_\Delta \,\cong\, 
  \Hom (E_{(\Delta, \partial\Delta)},A) \cr
\leftline{and}\cr
  E_{(\Delta, \partial \Delta)} 
  \;\buildrel \cong \over \longto\; 
  \D\E_{(\Delta, \partial\Delta)} 
  \,\cong\, \Hom (E_\Delta,A)\ .}
$$
If $\Delta$ is even quasi-convex, then the modules $E_\Delta$ and 
$E\b_{(\Delta, \partial \Delta)}$ are free, see [BBFK$_2$, 4.8, 
4.12].\qed

\medskip \noindent
Theorem 3.2 now follows from this proposition with $\F=\D\E$:

\bigskip
\noindent
{\bf 3.5 Proposition.} [BBFK${}_2$, 1.8 iii)] and [BreLu${}_2$, 3.14] 
{\it  For a fan $\Delta$ and two copies 
$\E\b$ and $\F\b$ of the equivariant intersection cohomology sheaf, 
every homomorphism 
$$
  \RR\b = E\b_o \,\to\, F\b_o = \RR\b \,
$$
extends in a unique manner to a homomorphism $\E\b \to \F\b$ of degree $0$.}

\medskip
\noindent
For its proof, we need a Vanishing Lemma. This is the place where the 
Hard Lefschetz Theorem enters:

\medskip
\noindent
{\bf 3.6 Lemma.} [BBFK${}_2$, 1.7, 1.8 ii)] and [BreLu${}_2$, 3.13] 
{\it For the 
equivariant intersection cohomology sheaf $\E\b$ on a 
non-zero cone $\sigma$, the following equivalent conditions hold:
\item{\rm(1)}\qquad
   $\overline E_\sigma^q = 0$ for $q \ge \dim \sigma$,
\item{\rm(2)}\qquad
   $\overline E_{(\sigma, \partial \sigma)}^q = 0$ for
   $q \le \dim \sigma$,
\item{\rm(3)}\qquad
   $E_{(\sigma, \partial \sigma)}^q = 0$ for $q \le \dim \sigma$.}

\medskip
\noindent
{\bf Proof:} We may assume $\dim\sigma = n$.\par
\noindent
(1) We use the setup of 0.3 (iii). First of all note that
$$
  B/\mm_B \otimes_B E_{\partial \sigma} \cong
  (B/\mm_B) [T] \otimes_{B [T]} E_{\partial \sigma}\ .
$$
Now we tensorize the exact sequence
$$
  0 \longto (B/\mm_B) [T] \buildrel \mu_T \over \longto (B/\mm_B) [T] 
  \longto A/\mm_A \longto 0
$$
with $E_{\partial \sigma}$ and obtain the exact sequence
$$
  (B/\mm_B) \otimes_B E_{\partial \sigma} 
  \buildrel \overline \mu_T \over \longto 
  (B/\mm_B) \otimes_B E_{\partial \sigma}
  \longto A/\mm_A \otimes_A E_{\partial \sigma} \longto 0
$$
with $\overline \mu_T:= \id_{(B/\mm_B)} \otimes \mu_T$, where $\mu_T$ 
acts on the $A$-module $E_{\partial \sigma}$. Thus
$$
  {\overline E}\b_{\partial \sigma} \;\cong\;
  \coker \Bigl(\overline \mu_T: 
    (B\b/\mm_{B\b}) \otimes_{B\b} E\b_{\partial\sigma}
    \longto
    (B\b/\mm_{B\b}) \otimes_{B\b} E\b_{\partial\sigma} 
  \Bigr)\,.
$$
On the other hand, according to [BBFK$_2$, (5.3.2)] together with 
(0.3.1) and using the notation of (1.3.1), we have an isomorphism
$$
  {\overline E}\b_{\partial \sigma} \;\cong\;
  \coker\Bigl(\overline \mu_{\psi}: IH\b (\Lambda_{\sigma}) \longto
  IH\b(\Lambda_{\sigma})\Bigr)\ ,
$$
where $\mu_{\psi}: \E (\Lambda_{\sigma}) \longto \E (\Lambda_{\sigma})$ 
is the multiplication with the strictly convex conewise linear function
$$
  \psi := T \,{\circ}\, (\pi|_{\partial \sigma})^{-1} 
  \;\in\; \A^2(\Lambda_{\sigma}) \,.
$$
It now suffices to apply for $m:=n{-1}$ the following theorem proved in [Ka]:

\medskip
\noindent
{\bf Hard Lefschetz Theorem.} {\it Let $\Lambda$ be a complete fan in 
the $m$-dimensional vector space $V$ and $\psi \in \A^2 (\Lambda)$ be a 
conewise linear strictly convex function. Then the homomorphism $L := 
\overline \mu_\psi$ induced by the multiplication $\mu_\psi: E_\Lambda 
\longto E_\Lambda$ with $\psi$ induces isomorphisms 
$$
  L^k: IH^{m-k}(\Lambda) \longto IH^{m+k}(\Lambda)
$$
for each $k \ge 0$. In particular $L$ is injective in degrees $q \le m-1$ 
and surjective in degrees $q \ge m-1$.}

\medskip
Let us finish the proof of 3.6: The equivalence of (1) and (2) follows 
from (2.6.1) and the dual pairing ($\quer {\rm PD}$) in Theorem 3.4 in 
the particular case $\Delta = \left<\sigma\right>$, while the 
equivalence of (2) and (3) is a consequence of this fact: For a finitely 
generated graded $A\b$-module $M\b$, one has $M^q=0$ for $q \le r$ if 
and only if $\overline M^q=0$ for $q \le r$.
\qed

\bigskip
\noindent
{\bf Proof of Proposition 3.5}:
For an inductive proof, we have to show that over each non-zero cone 
$\sigma$, a homomorphism $\phi_{\partial \sigma} \: E\b_{\partial \sigma}
\to F\b_{\partial \sigma}$ extends in a unique way to a homomorphism 
$\phi_{\sigma} \: E\b_{\sigma} \to F\b_{\sigma}$. By 3.6, (1), the 
$A\b$-modules $E\b_\sigma$ and $F\b_\sigma$ are generated by homogeneous 
elements of degree below $\dim \sigma$. On the other hand, 3.6,~(3) 
yields $E_{(\sigma, \partial \sigma)}^q = 0 = 
F_{(\sigma, \partial \sigma)}^q$ for $q \le \dim \sigma$. Hence, the 
restriction maps $E^q_\sigma \to E^q_{\partial \sigma}$ and $F^q_\sigma 
\to F^q_{\partial \sigma}$ are isomorphisms for $q < \dim\sigma$, whence 
the uniqueness of $\phi_\sigma$ follows. The existence is a consequence 
of the fact that $E_\sigma$ is a free $A_\sigma$-module. \qed

\bigskip
\bigskip\medskip\goodbreak
\centerline{\XIIbf 4. Comparison with previous definitions}

\medskip\nobreak
\noindent
Let  $\pi = \id_{V} : (V, \hat\Delta) \to (V, \Delta)$ be an {\it 
oriented refinement\/}, i.e., if a cone in $\hat \Delta^d$ is 
contained in a cone in $\Delta^d$, then their orientations coincide.

\bigskip
\noindent 
{\bf 4.1 Proposition.} {\it  For every pure sheaf $\F$ on 
$\hat \Delta$, there exists a canonical isomorphism} 
$$
  \D\bigl(\pi_*(\F)\bigr) \;\cong\; \pi_*(\D\F)\ .
$$ 

\medskip
\noindent 
{\bf Proof.}  For a cone  $\sigma \in \Delta^d$, let 
$\hat \sigma \preceq \hat \Delta$ denote its refinement. Then formula 
(2.3.1), applied to $\hat\sigma$, yields the isomorphism in the 
following chain 
$$
\eqalign{
  \D\bigl(\pi_*(\F)\bigr)_\sigma & = 
  \Hom (\pi_*(\F)_{(\sigma, \partial
  \sigma)}, A) \otimes \det V_\sigma^* =
  \Hom (F_{(\hat \sigma, \partial \hat
  \sigma)}, A) \otimes \det V_\sigma^*
   \cr
  & \cong 
  \D\F_{\hat \sigma} = \pi_*(\D\F)_\sigma\;.}\eqno{\qed}
$$

\medskip
 We now can prove the Compatibility Theorem:

\bigskip
\noindent
{\bf 4.2 Theorem.} [BreLu${}_2$, 7.2] {\it Let
$\hat\E\b$ be the equivariant intersection  cohomology sheaf of the oriented
refinement $\hat \Delta$ of the normal 
$n$-dimensional fan $\Delta$, and let $\iota : \E\b \to 
\pi_* (\hat \E\b)$, a homomorphism of  graded sheaves extending the 
identity $\E\b (o) = \RR = \pi_*(\hat \E\b)(o)$. Then the intersection  
products provide a commutative diagram}
$$
\matrix{
  \E\b (\Delta) \times 
  \E\b (\Delta, \partial\Delta)\kern-30pt &&
  \longto && 
  \kern-30pt\hat \E\b (\hat\Delta) \times 
  \hat \E\b (\hat\Delta,\partial \hat \Delta)\cr
  \noalign{\vskip6pt}
  & \searrow && \swarrow & \cr 
  \noalign{\vskip6pt}
  & & A\b[-2n]\,. &  &\cr} 
$$

\bigskip
\noindent
{\bf Proof.} The homomorphism $\iota$ provides a diagram 
$$
\matrix{
  \E
   &
  \buildrel \iota \over \longto &
  \pi_*(\hat \E)\cr
  \noalign{\vskip 6pt}
  \downarrow &
  & \downarrow \cr
  \noalign{\vskip 6pt}
   \D \E&
  \buildrel \D\iota \over \longleftarrow & 
  \D\pi_*(\hat \E) \cong \pi_*(\D\hat \E)
   \cr
  }
$$
where the vertical arrows are $\theta$ respectively $\pi_*(\hat \theta)$ 
with the duality correlations $\theta: \E \longto \D\E$ of (2.6.1) 
and $\hat \theta: 
\hat \E \longto \D \hat \E$. It is commutative at the zero cone and 
thus everywhere, see Proposition~3.5. Passing to the level of global 
sections yields the claim. \qed 

\bigskip
Finally let us discuss the approach of [BBFK$_2$, 6.1]. Here we use
the notion of an {\it evaluation map}:

\medskip
\noindent
{\bf 4.3 Definition.} {\it Let $\Delta$ be an oriented purely 
$n$-dimensional fan in the vector space $V$, endowed with a volume form 
$\omega \in \det V^*$. Then, for $1 \in E^0_\Delta = E^0_o =\RR$, the 
homomorphism 
$$
  e_\Delta^\omega := \chi_{\Delta}^\omega(1) : 
  E_{(\Delta, \partial \Delta)} \longto A[-2n] \;,
$$
see {\rm(3.3.1)}, is called the {\bfit evaluation map} associated 
to~$\omega$.}

\bigskip
\noindent
{\bf 4.4 Theorem.} {\it Let $\Delta$ be an oriented normal fan in a 
vector space $V$ endowed with a volume form $\omega \in \det V^*$. 
Furthermore let
$$
  \beta :\E\b \times \E\b \longto \E\b
$$
be a bilinear map of $\A$-module sheaves extending the multiplication 
$$
  E\b_o \times E\b_o = \RR \times \RR \longto \RR = E\b_o
$$
of real numbers. Then the pairing
$$
  e_\Delta^\omega \circ \beta_{\Delta} :
  \E (\Delta) \times \E (\Delta, \partial \Delta)
  \longto \E (\Delta, \partial \Delta) \longto A[-2n] 
  \leqno(4.4.1)
$$
coincides with the intersection  product.}

\bigskip
\noindent
Note that for a simplicial fan $\Delta$, the equality $\E\b = \A\b$ 
holds, so the bilinear map $\beta$ necessarily is the multiplication of 
functions and thus, symmetric. In the non-simplicial case, the 
map~$\beta$ is not uniquely determined. Nevertheless, there always 
exists such a map~$\beta$ that is symmetric. For a complete fan, the 
intersection product is thus symmetric, which also follows from Theorem 
4.2 with a simplicial subdivision $\hat \Delta$ of~$\Delta$.

\bigskip
\noindent
{\bf Proof.} For each cone $\sigma$ and a positive volume form 
$\omega_\sigma \in \det (V_\sigma^*)$, we define 
$e_{\sigma}^{\omega_{\sigma}}$ analogously to $e_\Delta^\omega$. Then 
$e_\sigma^{\omega_\sigma} \otimes \omega_\sigma$ does not depend on the 
choice of $\omega_\sigma$, and the family of homomorphisms  
$$
  \tilde\phi_{\sigma} : \E_{\sigma} \longto \D\E_{\sigma} 
  \;,\quad
  s \longmapsto 
  \bigl (e_{\sigma}^{\omega_{\sigma}} \circ \beta \bigr )
   (s, \_ ) \otimes\, \omega_{\sigma} \;
$$
defines a duality correlation $\tilde\phi : \E\b \to \D\E\b$ and thus, 
according to Theorem 3.2, is unique. In particular the pairing (4.4.1) 
is the intersection product. \qed

%
\bigskip
{\centerline{\bf References}}
\medskip
\def\liitem{\par\noindent
                   \hangindent=40pt\ltextindent}
\def\ltextindent#1{\hbox to \hangindent{#1\hss}\ignorespaces}

\liitem{[BBFK$_{1}$]} \quad{\sc G.~Barthel, J.-P.~Brasselet,
K.-H.~Fieseler, L.~Kaup:} {\it Equivariant
Intersection Cohomology of Toric Varieties}, Algebraic
Geometry, Hirzebruch 70, 45--68, Contemp. Math. {\bf 241},
Amer. Math. Soc., Providence, R.I., 1999.

\liitem{[BBFK$_{2}$]} {\sc --, --, --, --:} 
{\it Combinatorial Intersection Cohomology 
for Fans}, T\^ohoku Math.~J. {\bf 54} (2002), 
1--41.

\liitem{[BreLu$_1$]} {\sc P.~Bressler, V.~Lunts:} {\it
Intersection cohomology on nonrational polytopes},
(pr)e-print {\tt math.AG/0002006} (33~pages), 2000,
to appear in Compositio Math.

\liitem{[BreLu$_2$]} {\sc --, --:} {\it Hard Lefschetz theorem and 
Hodge-Riemann relations for intersection cohomology of nonrational 
polytopes}, (pr)e-print {\tt math.AG/0302236~v2} (46~pages), 2003.

\liitem{[Bri]} {\sc M.~Brion:} {\it The structure of the polytope 
algebra}, T\^ohoku Math.~J. {\bf 49} (1997), 1--32.

\liitem{[Fi]} {\sc K.-H.~Fieseler:} {\it Combinatorial Duality and
Intersection Product}, (pr)e-print {\tt math.AG/0306344 v2} (12~pages), 
2003.

\liitem{[Ka]} {\sc K.~Karu:} {\it Hard Lefschetz Theorem for Nonrational 
Polytopes}, (pr)e-print {\tt math.AG/0112087} (25~pages), 2002.

\liitem{[Mc]} {\sc P.~McMullen:} {\it On simple Polytopes}, Invent.~math. 
{\bf 113} (1993), 419-444.

\liitem{[St]} {\sc R.~Stanley} {\it Generalized h-vectors, intersection 
cohomology of toric varieties and related results}, in {\sc M.~Nagata, 
H.~Matsumura} (eds.): {\it Commutative Algebra and Combinatorics},  
Adv.\ Stud.\ Pure Math.~{\bf 11}, Kinokunia, Tokyo, and North Holland, 
Amsterdam/New York, 1987, 187-213.

\bigskip\bigskip\goodbreak\noindent

\vbox{
{\klein Addresses of authors}
\bigskip
\nobreak
\hbox to \hsize{\hfill
{\klein
\parskip 0pt
\baselineskip=12pt
\vtop {\vskip 0.5truecm} {\hskip 0.1truecm}
\vbox{\hsize=0.5\hsize
\obeylines{G. Barthel, L. Kaup
Fachbereich Mathematik
\quad und Statistik
Universit\"at Konstanz
Fach D 203
D-78457 Konstanz
e-mail:
{\tt Gottfried.Barthel@uni-konstanz.de
Ludger.Kaup@uni-konstanz.de}
}}\hfill
\vbox{\hsize=0.45\hsize
\obeylines{K.-H. Fieseler
Matematiska Institutionen
Box 480
Uppsala Universitet
SE-75106 Uppsala
e-mail: {\tt khf@math.uu.se}
{\ }
J.P. Brasselet.
IML/CNRS, Luminy Case 907
F-13288 Marseille Cedex 9
e-mail: {\tt jpb@iml.univ-mrs.fr}}}\hfill}}
}

\bye

\bye